\documentclass[12pt]{amsart}
\usepackage{amsfonts,amsthm,amsmath,amssymb,verbatim}
\usepackage{graphicx}
\usepackage[active]{srcltx}
\usepackage[all,cmtip]{xy}

\newtheorem{thm}{Theorem}[section]
\newtheorem{prop}[thm]{Proposition}

\newtheorem{conj}[thm]{Conjecture}

\theoremstyle{remark}

\newtheorem{example}[thm]{Example}
\newtheorem{remark}[thm]{Remark}
\newtheorem{defin}{Definition}



\def\C{\mathbb{C}}
\def\R{\mathbb{R}}

\def\Z{\mathbb{Z}}
\def\P{\mathbb{P}}

\def\Sch{\mathfrak{S}}

\def\Sym{{\rm Sym}}

\def\Volume{{\rm Volume}}

\def\wt{\widetilde}

\def\l{\lambda}
\def\G{\Gamma}

\def\emptyset{\varnothing}

\title{Schubert calculus on Newton--Okounkov polytopes}
\author{Valentina Kiritchenko}
\email{vkiritch@hse.ru}
\thanks{The study has been partially funded by the Russian Academic Excellence Project '5-100'.}

\address{Laboratory of Algebraic Geometry and Faculty of Mathematics\\
National Research University Higher School of Economics\\
Usacheva str. 6, 119048 Moscow, Russia}
\address{Institute for Information Transmission Problems, Moscow, Russia}
\author{Maria Padalko}
\date{}

\begin{document}
\begin{abstract}
A Newton--Okounkov polytope of a complete flag variety can be turned into a convex geometric model for Schubert calculus.
Namely, we can represent Schubert cycles by linear combinations of faces of the polytope so that the intersection product of cycles corresponds to the set-theoretic intersection of faces (whenever the latter are transverse).
We explain the general framework and survey particular realizations of this approach in types $A$, $B$ and $C$.
\end{abstract}

\maketitle

\section{Introduction}
Theory of Newton--Okounkov convex bodies \cite{KK,LM} allows us to
apply ideas of toric geometry in the non-toric setting.
In this paper, we explore non-toric applications of {\em polytope
rings} (see Section \ref{s.prelim} for a definition) introduced by Khovanskii and Pukhlikov \cite{PKh}.
With a convex polytope $P\subset\R^d$, they associated a graded commutative ring (the polytope ring):
$$R_P=\bigoplus_{i=0}^d R_P^i$$
that has Poincar\'e duality.
The polytope rings were originally used to give a convenient functorial description of the cohomology rings of smooth toric varieties.
In this case, $P$ is always a simple lattice polytope, that is, all vertices of $P$ belong to $\Z^d\subset\R^d$, and only $d$ edges meet at every vertex of $P$.
In \cite{Kaveh}, Kaveh noted that polytope rings can also be used for a partial description of the cohomology rings of
spherical varieties.
In this case, $P$ is still a lattice polytope but not necessarily simple.

For simple polytopes, every face $\G\subset P$ can be naturally identified with an element $x_\G\in R_P$ so that $$x_\G x_{\G'}=x_{\G\cap\G'}$$
for any two transverse faces $\G$ and $\G'$.
This is no longer true for non-simple polytopes, that is, individual faces of $P$ do not have natural counterparts in $R_P$.
However, it is still possible to identify every element of $R_P$ with a linear combination of faces of $P$ so that the product in the polytope ring corresponds to the intersection of faces.
In \cite{KST}, the first author, Smirnov and Timorin developed a general framework for such calculus on polytopes, and studied its applications to Schubert calculus on Gelfand--Zetlin polytopes in type $A$.
In this paper, we mainly consider applications to Schubert calculus in types $B$ and $C$.

Representation theory of classical  groups is a source of several interesting families
of lattice convex polytopes.
For $SL_n(\C)$ (type $A$), there is a well-known family of Gelfand--Zetlin (GZ) polytopes $GZ_\lambda$.
Here
$\l:=(\l_1,\ldots,\l_n)\in\Z^n$ runs through dominant weights of $SL_n(\C)$, that is, $\l_1\ge\l_2\ge\ldots\ge\l_n$.
Originally, GZ polytopes were constructed using representation theory, namely, lattice points in the polytope
$GZ_\lambda$ parameterize the vectors in a special basis in the irreducible representation $V_\l$ of $SL_n(\C)$
with the highest weight $\l$ (see \cite{M} for a survey on GZ bases).
In convex geometric terms, the GZ polytope $GZ_\l\subset\R^d$, where $d:=\frac{n(n-1)}{2}$, is defined as the set of all points
$(z^1_1,z^1_2,\ldots, z^1_{n-1};z^2_1,\ldots,z^2_{n-2};\ldots; z^{n-1}_1)\in\R^d$ that satisfy the following
interlacing inequalities:
$$
\begin{array}{cccccccccc}
\l_1&       & \l_2    &         &\l_3          &    &\ldots    & &       &\l_n   \\
    &z^1_1&         &z^1_2  &         & \ldots   &       &  &z^1_{n-1}&       \\
    &       &z^2_1 &       &  \ldots &   &        &z^2_{n-2}&         &       \\
    &       &       &  \ddots   & &  \ddots   &      &         &         &       \\
    &       &       &  &z^{n-2}_1&     &  z^{n-2}_2 &        &         &       \\
    &       &         &    &     &z^{n-1}_1&   &              &         &       \\
\end{array}
\eqno(GZ_A)$$
where the notation
$$
 \begin{array}{ccc}
  a &  &b \\
   & c &
 \end{array}
 $$
means $a\ge c\ge b$ (the table encodes $2d$ inequalities).
Figure 1 shows the $3$-dimensional GZ polytope for $n=3$ and $\l=(3,0,-3)$.
Note that $GZ$ polytopes are not simple.

\begin{figure}
\begin{center}
\includegraphics[width=7cm]{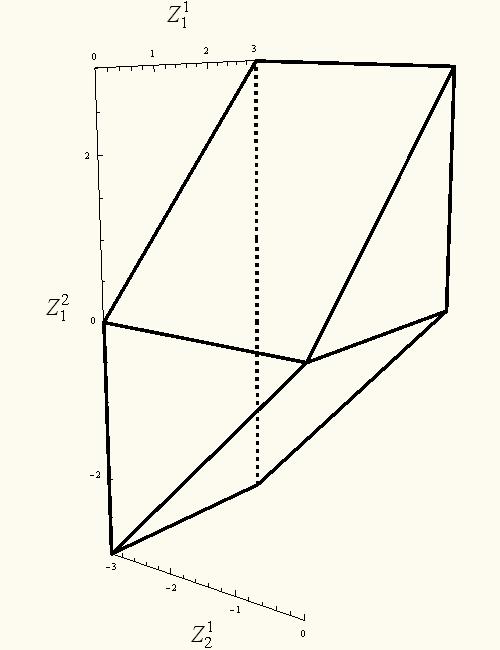}
\caption{GZ polytope in type $A$ for $n=3$ and $\l=(3,0,-3)$}
\end{center}
\end{figure}

GZ polytopes in types $B$, $C$ and $D$ were defined in
\cite{BZ} (see Section \ref{ss.GZ} for definitions in types $B$ and $C$) and
are related to representation theory of $SO_{2n+1}(\C)$, $Sp_{2n}$ and $SO_{2n}(\C)$, respectively.
They are special cases of {\em string polytopes} introduced by Berenstein--Zelevinsky and
Littelmann \cite{L}.
There are other families of polytopes in representation theory such as Nakashima--Zelevinsky
polyhedral realizations of crystal bases and Feigin--Fourier--Littelmann--Vinberg polytopes.
They have representation-theoretic meaning similar to that of string polytopes
but are not combinatorially equivalent to the latter.
All these polytopes were exhibited as Newton--Okounkov polytopes of complete flag varieties for certain
geometric valuations \cite{FFL14,FFL15,FO,Ka,K15} (see Section \ref{ss.NO} for more details).

For $G=SL_n(\C)$, the complete flag variety $G/B$ (here $B\subset G$ denotes the subgroup of upper-triangular matrices) can be thought of as a variety of complete flags of subspaces $(\{0\}\subset V^1 \subset V^2\subset\ldots\subset V^{n-1}\subset\C^n)$
where $\dim V^i=i$, and there are no gaps.
There are similar descriptions of complete flag varieties $G/B$ for other classical groups $G$ (see Section \ref{ss.NO}).
Recall that globally generated line bundles $L_\l$ on $G/B$ are in bijective correspondence with irreducible representations $V_\l$ of $G$
so that $H^0(L_\l, G/B)\simeq V_\l^*$ \cite[Proposition 1.4.5]{Brion}.
Here $\l$ runs through the dominant weights of $G$.
We denote by $\deg_\l (G/B)$ the degree of the image of $G/B$ under the map $G/B\to \P(V_\l)=\P(H^0(L_\l, G/B)^*)$.

In \cite{Kaveh}, polytope rings of string polytopes were identified with the cohomology rings of complete flag varieties.
More generally, string polytope in this description can be replaced with any linear family (in the sense of \cite{KaVi}) of convex polytopes $P_\l$ parameterized by the dominant weights $\l$ whenever the following identity holds:
$$\Volume(P_\l)=d!\deg_\l(G/B)\eqno(1)$$
where $d:=\dim G/B$.
We regard both sides of this identity as polynomials in $\l$.
In particular, polytopes $P_\l$ yield an analog of Kushnirenko's theorem for $G/B$.

Since Newton--Okounkov polytopes of line bundles on $G/B$ by construction satisfy identity (1)
they can be used to model Schubert calculus.
Recall that the cohomology ring $H^*(G/B,\Z)$ has a special basis of {\em Schubert cycles} $[X_w]$ with striking positivity properties.
Namely, the structure constants (i.e., the coefficients $c^u_{v w}$ in the decomposition $[X_w][X_v]=\sum_u c^u_{v w}[X_u]$) are always non-negative.
However, no enumerative meaning (in the spirit of Littlewood--Richardson rule for Grassmannians) of these coefficients is known.
Polytope rings provide a new framework for combinatorial interpretation of structure constants.
An important task is to find presentations of Schubert cycles in polytope rings by linear combinations of faces with positive coefficients.
Another task is to find Newton--Okounkov polytopes for which these presentations have especially simple combinatorics.
It is tempting to use Grossberg--Karshon cubes \cite{HY16,HY17} since they are combinatorial cubes.
However, there are might be issues with positivity, that is, some Schubert cycles will be represented by linear combinations
of faces with negative coefficients (see Example \ref{e.GK}).

There is an algorithm (geometric mitosis) for finding positive presentations of Schubert cycles by faces using convex geometric analogs of Demazure operators from representation theory \cite{K16I,K16II}.
In the present paper, we describe geometric mitosis in more combinatorial terms, outline its applications and formulate conjectures.
For GZ polytopes in type $A$, this algorithm reduces to Knutson--Miller mitosis on pipe dreams and was used in \cite{KST}.
In types $B$ and $C$, geometric mitosis reduces to a different combinatorial rule that conjecturally yields presentations of Schubert cycles by faces of GZ polytopes in respective types.
In particular, 4-dimensional GZ polytope in type $C_2$ can be used to model Schubert calculus on the variety of isotropic flags
in $\C^4$ \cite{P17}.
Another convex geometric model for the same flag variety was constructed in \cite{I} using a different string polytope in type $C_2$.

\section{Preliminaries}\label{s.prelim}
In this section, we recall the definitions of polytope rings, GZ polytopes and flag varieties in types $B$ and $C$.
We discuss the relationship between the polytope rings of GZ polytopes and cohomology rings of flag varieties.
We also define Newton--Okounkov polytopes of flag varieties.

\subsection{Polytope ring}
Let $L\subset \R^d$ be a lattice, and $P\subset \R^d$ a convex polytope whose vertices lie in $L$.
We say that $P$ is a {\em lattice polytope} with respect to $L$.
By the {\em standard lattice} $\Z^d$ we mean the lattice
$\{(x_1,\ldots,x_d)\in\R^d\ |\ x_i\in \Z \mbox{ for all } i=1,\ldots, d \}$.
We choose the translation invariant volume form on $\R^d$ so that the covolume of $L$ is $1$.

Recall that two convex polytopes $P$ and $Q$ are called {\em analogous} if they have the same normal fan,
i.e. there is a one-to-one correspondence between the faces of $P$ and the
faces of $Q$ such that any linear functional, whose restriction to $P$
attains its maximal value at a given face $F\subseteq P$ has the
property that its restriction to $Q$ attains its maximal value at the corresponding face of $Q$.

Denote by $S_P$ the set of all polytopes analogous to $P$.
This set can be endowed with the structure of a commutative semigroup using {\em Minkowski sum}
$$
P_1+P_2=\{x_1+x_2\in\R^d\ |\ x_1\in P_1,\ x_2\in P_2\}
$$
It is not hard to check that this semigroup has cancelation property.
We can also multiply polytopes in $S_P$ by positive real numbers using dilation:
$$
\lambda P=\{\lambda x\ |\ x\in P\},\quad \lambda\ge 0.
$$
Hence, we can embed the semigroup of convex polytopes into its Grothendieck group $V_P$,
which is a real vector space.
The elements of $V_P$ are called {\em virtual polytopes} analogous to $P$.

On the vector space $V_P$, there is a homogeneous polynomial $vol_P$ of degree $d$, called the
{\em volume polynomial}.
It is uniquely characterized by the property that its
value $vol_P(Q)$ on any convex polytope $Q\in S_P$ is equal to the volume of $Q$.

Let $\Lambda_P$ be a lattice in $V_P$ generated by some lattice polytopes (with respect to $L$)
analogous to $P$ (we do not assume that $\Lambda_P$ contains all
lattice polytopes analogous to $P$).
The symmetric algebra $\Sym(\Lambda_P)$ of $\Lambda_P$ can be
thought of as the ring of differential operators with constant
integer coefficients acting on $\R[V_P]$, the space of all polynomials on $V_P$.
If $D\in\Sym(\Lambda_P)$ and $\varphi\in\R[V_P]$, then we write
$D\varphi\in\R[V_P]$ for the result of this action.
Define $A_P$ as the homogeneous ideal in $\Sym(\Lambda_P)$ consisting of all
differential operators $D$ such that $Dvol_P=0$.
Set $R_P=\Sym(\Lambda_P)/A_P$.
This ring is called {\em the polytope ring} associated with the polytope $P$ and the lattice
$\Lambda_P$.

\begin{example}
Let $L=\Z^d$ be the standard lattice, and $P$ an integrally simple lattice polytope (that is, only $d$ edges meet at every vertex of $P$, and primitive vectors on these edges span $L$ over $\Z$).
Let $\Lambda_P$ be the lattice in $V_P$ generated by all lattice polytopes (with respect to $L$)
analogous to $P$.
Then the ring $R_P$ is isomorphic to the Chow (or cohomology) ring $H^*(X_P,\Z)$ of the smooth toric variety $X_P$ associated with the normal fan of $P$ \cite{PKh}.

When $P$ is simple, every facet $\Gamma\subset P$ defines a differential operator $\partial_\Gamma\in R_P$ (see \cite[Section 2.3]{KST} for the details).
Recall that the closures of torus orbits in $X_P$ are in bijective correspondence with faces of $P$.
They also give a generating set in the cohomology ring $H^*(X_P,\Z)$.
Every face $F=\Gamma_1\cap\ldots\cap\Gamma_k$ can be identified with the operator $[F]=\partial_{\Gamma_1}\cdots\partial_{\Gamma_k}\in R_P$.
Using linear relations between $\partial_\Gamma$ in $R_P$ we can compute products in $H^*(X_P,\Z)$ by intersecting faces of $P$.

For instance, if $P\subset\R^2$ is the trapezoid with vertices $(0,0)$, $(1,0)$, $(0,1)$, $(1,2)$, then the corresponding toric variety $X_P$ is the blow-up of $\C\P^2$ at one point.
The edge $\G_1=\{x=0\}$ corresponds to the exceptional divisor $E\subset X_P$.
The other edges are $\G_2=\{y-x=1\}$, $\G_3=\{ x=1\}$ and $\G_4=\{y=0\}$.
There are two linear relations between $\partial_{\Gamma_i}$.
Namely, the parallel translations along $x$ and $y$ axes do not change the area of $P$, hence,
$\partial_{\Gamma_1}+\partial_{\Gamma_2}=\partial_{\Gamma_3}$ and $\partial_{\Gamma_2}=\partial_{\Gamma_4}$.
In particular, the identity $[E]^2=-[pt]$ in $H^*(X_P,\Z)$ can be obtained as follows:
$$[\G_1]^2=[\G_1]([\G_3]-[\G_2])=[\G_1\cap\G_3]-[\G_1\cap\G_2]=-[pt].$$

\end{example}

\begin{example} Let $L=\Z^d$, and $P$ the GZ polytope in type $A$ corresponding to a strictly dominant $\l=(\l_1,\ldots,\l_n)$ (that is, $\l_1>\l_1>\ldots>\l_n$).
Let $\Lambda_P$ be the lattice in $V_P$ generated by all GZ polytopes $P_\l$ for all dominant $\l$.
Then the ring $R_{GZ}:=R_P$ is isomorphic to the cohomology ring $H^*(GL_n(\C)/B,\Z)$ of the complete flag variety in type $A$ \cite{Kaveh}.

Since the GZ polytope is not simple, there is no correspondence between individual faces of $P$
and elements of $H^*(GL_n(\C)/B,\Z)$.
However, it is possible to identify every element of
$H^*(GL_n(\C)/B,\Z)$ with a linear combination of faces of $P$ (see \cite[Section 2]{KST} for more details).
Again, we can compute all products in $H^*(GL_n(\C)/B,\Z)$ by intersecting faces of $P$ (see
\cite[Section 2.4]{KST} for an example of such computations).
\end{example}

In what follows, $L$ will be a sublattice of
$\frac{1}{2}\Z^d:=\{(x_1,\ldots,x_d)\ |\ 2x_i\in \Z \mbox{ for all } i=1,\ldots, d \}$.
We always compute volumes of faces of $P$ with respect to the lattice $L$.
More precisely, if
$F\subset P$ is a face, and $\R F$ is its affine span then the volume of the face is computed using the volume form on $\R F$ normalized so that the covolume of $L\cap \R F$ is $1$.

\subsection{GZ polytopes in types $B$ and $C$}\label{ss.GZ}
Let $\lambda=(\lambda_1,\ldots,\lambda_n)$ be a non-increasing collection of non-negative integers.
Put $d=n^2$.
Denote coordinates in $\R^d$ by $(x^1_1,\ldots, x^1_n; y^1_1,\ldots, y^1_{n-1};\ldots; x_1^{n-1}, x^{n-1}_2, y^{n-1}_1; x^n_1)$.
For every $\l$, define the {\em symplectic GZ polytope} $SGZ_\l\subset\R^d$ for $Sp_{2n}(\C)$ by the following interlacing inequalities:
$$
\begin{array}{ccccccccccc}
\l_1&       & \l_2  &      &\l_3   &       & \ldots   & \l_n    &         &0      &\\
    &x^1_1  &       &x^1_2 &       &\ldots &          &         &x^1_{n}  &       &\\
    &       &y^1_1  &      &y^1_2  &       & \ldots   &y^1_{n-1}&         & 0     &\\
    &       &       &x^2_1 &       &\ldots &          &         &x^2_{n-1}&       &\\
    &       &       &      & y^2_1 &       & \ldots   &y^2_{n-2}&         & 0     &\\
    &       &       &      &       &\ddots &          & \vdots  &         &\vdots &\\
    &       &       &      &       &       & x^{n-1}_1&         &x^{n-1}_2&       &\\
    &       &       &      &       &       &          &y^{n-1}_1&         &     0 &\\
    &       &       &      &       &       &          &         & x^{n}_1 &       &\\
\end{array}
\eqno{GZ_C}$$
Again, every coordinate in this table is bounded from above by its upper left neighbor and bounded from below by its upper right neighbor (the table encodes $2d$ inequalities).
We regard $SGZ_\l$ as a lattice polytope with respect to the standard lattice $\Z^d$.
Roughly speaking, $SGZ_\l$ is the polytope defined using half of the GZ pattern $(GZ_A)$ for $SL_{2n}(\C)$.
\begin{example}\label{e.SGZ}
The polytope $SGZ_\l\subset\R^4$ for $Sp_4(\C)$ is given by 8 inequalities:
$$\l_1\ge x^1_1\ge \l_2; \quad \l_2\ge x^1_2\ge 0; \quad  x^1_1\ge y^1_1\ge x^1_2; \quad y^1_1\ge x^2_1\ge 0.$$
It is not hard to compute the volume polynomial of $SGZ_\l$:
$$vol_{SGZ}(\l_1,\l_2)=\frac16\l_1\l_2(\l_1-\l_2)(\l_1+\l_2).$$
This volume times $4!$ is equal to the degree $\deg_\l (Sp_4(\C)/B)$ of the isotropic flag variety.
\end{example}

The polytope ring $R_{SGZ}$ defined by the family of symplectic GZ polytopes is isomorphic to the cohomology ring $H^*(Sp_{2n}(\C)/B,\Z)$.
Indeed, by \cite{Kaveh} it is isomorphic to the subring of $H^*(Sp_{2n}(\C)/B,\Z)$ generated by the first Chern classes of line bundles $L_\l$ corresponding to the weights of $Sp_{2n}(\C)$.
Since the torsion index of $Sp_{2n}(\C)$ is $1$, this subring coincides with the whole ring (see \cite{T} for the details on torsion indices of classical groups).

The {\em odd orthogonal GZ polytope} $OGZ_\l\subset\R^d$ for $SO_{2n+1}(\C)$ is defined using the same pattern $(GZ_C)$ but a different lattice $L_B\subset\R^d$.
Namely, $L_B$ consists of all points $(x^1_1,\ldots, x^1_n; y^1_1,\ldots, y^1_{n-1};\ldots; x_1^{n-1}, x^{n-1}_2, y^{n-1}_1; x^n_1)\in\frac12\Z^d$ such that all coordinates except for $x^1_n$, $x^2_{n-1}$,\ldots, $x^n_1$ are integer.
Lattice points $SGZ_\l\cap\Z^d$ and $SGZ_\l\cap L_B$ parameterize basis vectors in irreducible representations of $Sp_{2n}(\C)$
and $SO_{2n+1}(\C)$, respectively (see \cite[Section 6]{L} for more details).

\begin{remark} Family of odd orthogonal GZ polytopes (as defined in \cite{BZ,L}) consists of two subfamilies parameterized by integer and half-integer $\l$.
The group $SO_{2n+1}(\C)$ is not simply connected, and half-integer weights correspond to the characters of the maximal torus in the universal cover $Spin(2n+1)$.
If we define the polytope ring $R_{SGZ}$ using  the first subfamily we get a subring of $H^*(SO_{2n+1}/B,\Z)$ generated by the first Chern classes of line bundles $L_\l$ corresponding to the characters $\l$ of the maximal torus in $SO_{2n+1}(\C)$.
\end{remark}

\begin{example}\label{e.B_2}
The polytope $OGZ_\l\subset\R^4$ for $Sp_4(\C)$ is given by the same 8 inequalities as in Example \ref{e.SGZ}.
However, its volume polynomial is computed using a different volume form chosen so that the covolume of $L_B$ is 1.
Since $\Z^4\subset L_B$ has index $4$, we get $vol_{OGZ}=4vol_{SGZ}$.

There is an exceptional isomorphism $Sp_4(\C)/{\pm 1}\simeq SO_5(\C)$.
In particular, flag varieties in types $B_2$ and $C_2$ are the same.
This isomorphism takes the dominant weight $\l=(\l_1,\l_2)$ of $Sp_4(\C)$ to the dominant
weight $\wt \l=(\l_1+\l_2)/2,(\l_1-\l_2)/2)$ of $SO_5(\C)$.
This agrees with the identity $vol(SGZ_\l)=vol(OGZ_{\wt \l})$.
\end{example}

\subsection{Newton--Okounkov polytopes of flag varieties}\label{ss.NO}
We recall a definition of Newton--Okounkov convex bodies in the case of flag varieties.
We refer the reader to \cite{KK,LM} for definitions in the more general setting.

Recall that the complete flag variety $SL_n(\C)/B$ is defined as the variety of complete flags of subspaces $M^\bullet=(\{0\}\subset V^1 \subset V^2\subset\ldots\subset V^{n-1}\subset\C^n)$.
We define $SO_n(\C)/B$ and $Sp_{2n}/B$ as subvarieties of {\em orthogonal} and {\em isotropic} flags in $SL_n(\C)/B$ and $SL_{2n}/B$, respectively.
A complete flag $M^\bullet$ in $\C^n$ is {\em orthogonal} if $V^i$ is orthogonal to to $V^{n-i}$ with respect to a non-degenerate symmetric bilinear form fixed by $SO_n(\C)$.
Let $\omega$ be a non-degenerate skew-symmetric bilinear form fixed by $Sp_{2n}(\C)$.
A complete flag $M^\bullet$ in $\C^{2n}$ is called {\em isotropic} if the restriction of $\omega$ to $V^n$ is zero, and $V^{2n-i}=\{v\in \C^{2n} \ | \ \omega(v,u)=0 \mbox{ for all } u\in V^i\}$.

Every flag variety $X$ of dimension $d$ has an open dense subset $C$ ({\em open Schubert cell}) isomorphic to the affine space $\C^d$.
It can be constructed as follows.
Fix a complete flag
$F^\bullet:=(F^1\subset F^2\subset\ldots\subset F^{n-1}\subset \C^n)$ such that $F^\bullet\in X$ (this amounts to fixing a Borel subgroup $B\subset G$).
Also fix a basis $e_1$,\ldots, $e_n$ in $\C^n$ compatible with $F^\bullet$ (or a maximal torus in $B$), that is, $F^i=\langle e_1,\ldots,  e_i\rangle$.
{\em The open Schubert cell} $C$ with respect to $F^\bullet$ is defined as the set of all flags $M^\bullet$ that are in general position with the standard flag $F^\bullet$, i.e., all  intersections $M^i\cap F^j$ are transverse.
Let $x_1$, \ldots, $x_d$ be coordinates on the open Schubert cell $C$.

\begin{example}\label{e.NO_A} In type $A$, we can identify the open Schubert cell  $C$  with an affine space $\C^d$ (for $d=\frac{n(n-1)}{2}$) by choosing for every flag $M^\bullet$  a basis  $v_1$,\ldots, $v_n$ in $\C^n$ of the form:
$$v_1=e_n+x^{n-1}_1 e_{n-1}+\ldots+ x^1_1 e_1, $$
$$v_2=e_{n-1}+x^{n-2}_2 e_{n-2}+\ldots+ x^1_2 e_1, \quad \ldots \quad, v_{n-1}=e_2+x^1_{n-1}e_1, \quad v_n=e_n,$$
so that $M^i=\langle v_1,\ldots,  v_i\rangle$.
Such a basis is unique, hence, the coefficients $(x^i_j)_{i+j< n}$ are coordinates on the open cell.
In other words, every flag $M^\bullet\in C$ gets identified with a triangular matrix:
$$\begin{pmatrix}
   x^1_1 & x^1_2& \ldots & x^1_{n-1} & 1\\
   x^2_1 & x^2_2   & \ldots & 1 & 0\\
   \vdots & \vdots   &  &  & 0\\
   x^{n-1}_1 & 1 & \ldots & 0 & 0\\
    1        & 0 &        & 0 & 0\\
  \end{pmatrix} \eqno(FFLV).
$$
Similar coordinates can be introduced on flag varieties in other types.
\end{example}

Let $V\subset \C(X)=\C(x_1,\ldots,x_d)$ be a finite-dimensional subspace of rational functions on $X$.
Our main examples are spaces of global sections $H^0(L_\l, X)\simeq V_\l^*$ of line bundles on $X$.
We fix a section $s_0\in H^0(L_\l,X)$, and identify sections $s\in H^0(L_\l,X)$ with rational functions $f=\frac{s}{s_0}\in\C(X)$.

\begin{example} We continue Example \ref{e.NO_A}.
If
$$\l=(\underbrace{1,\ldots,1}_k,\underbrace{0\ldots,0}_{n-k}),$$ then $V_\l^*$ can be identified with the subspace of $\C(x^i_j)_{i+j< n}$ spanned by the minors of the $n\times k$ matrix formed by the first $k$ columns of the matrix $(FFLV)$.
These minors are exactly the Pl\"ucker coordinates of the Grassmannian $G(k,n)$ in the Pl\"ucker embedding.
The map $X\to H^0(L_\l,X)^*$ is the composition of the projection $X\to G(k,n)$ (obtained by forgetting all subspaces in the flag $M^\bullet$ except for the $V^k$) and the Pl\"ucker embedding of $G(k,n)$.
\end{example}

To assign the {\em Newton--Okounkov convex body} to $V$ we need an extra ingredient.
Choose a translation-invariant total order on the lattice $\Z^d$ (e.g., we can take the lexicographic order).
Consider a map
$$v:\C(x_1,\ldots,x_d)\setminus\{0\}\to\Z^d,$$
that behaves like the lowest order term of a polynomial, namely: $v(f+g)\ge \min\{v(f),v(g)\}$ and $v(fg)=v(f)+v(g)$ for all nonzero $f,g$.
Recall that maps with such properties are called {\em valuations}.

\begin{defin}
The {\em Newton--Okounkov convex body} $\Delta_v(X,V)$ is the closure
of the convex hull of the set
$$\bigcup_{k=1}^\infty\left\{\frac{v(f)}{k} \ |\ f\in V^k \right\}\subset\R^d.$$
By $V^k$ we denote the subspace spanned by the $k$-th powers of the functions from $V$.
\end{defin}

\begin{example}
Using coordinates of Example \ref{e.NO_A} we can define the valuation $v$ as follows.
Order the coefficients $(x^i_j)_{i+j< n}$ of the matrix $(FFLV)$ by starting from column $(n-1)$ and going from top to bottom in every column and from right to left along columns.
Then $\Delta_v(X,V_\l^*)$ coincides with the Feigin--Fourier--Littelmann--Vinberg polytope $FFLV(\l)$ \cite{K15}.
Moreover, the inclusion $FFLV(\l)\subset\Delta_v(X,V_\l^*)$ follows from a straightforward computation of the valuation $v$ on the minors of the matrix $(FFLV)$ (see \cite[Example 2.9]{K15} for more details).
\end{example}

Different valuations might yield different Newton--Okounkov convex bodies.
In particular, GZ polytopes can also be obtained as Newton--Okounkov polytopes of flag varieties \cite{Ka,FO}.
Okounkov made the first explicit computation of this kind, namely, he exhibited symplectic GZ polytopes as Newton--Okounkov polytopes of the isotropic flag varieties \cite{O}.

\section{Geometric mitosis}\label{s.mitosis}
In \cite{K16II}, convex geometric analogs of Demazure (or divided difference) operators are defined on convex polytopes and used to construct {\em DDO polytopes} that have the same properties as Newton--Okounkov polytopes of flag varieties.
In \cite{K16I}, operations on faces of a DDO polytope (geometric mitosis) are defined that yield positive presentations of Schubert cycles by faces.
Here we define the same operations in more combinatorial terms using a vertex cone instead of a DDO polytope.
We refer the reader to \cite[Theorem 3.6]{K16II}, \cite[Proposition 2.5]{K16I} for connections with representation theory and convex geometry.

\begin{example}\label{e.GK}
Figure 2 illustrates the idea of mitosis in the simplest example.
The trapezoid and rectangle on the left picture have the same number of lattice points with given sum of coordinates.
The same is true for the right picture.
However, the trapezoid on the right picture becomes a virtual polytope (in particular, lattice points marked with circles has to be counted with the zero coefficient) while the rectangle remains a true polytope.
There is a price to pay: the left vertical edge of the trapezoid corresponds to two edges of the rectangle (that is, a single edge of the trapezoid has the same number of lattice points as two edges of the rectangle).
In short, mitosis preserves positivity at the cost of more involved combinatorics.
\end{example}

Consider a vector space with the direct sum decomposition
$$\R^d=\R^{d_1}\oplus\ldots\oplus\R^{d_r},$$
and choose coordinates $x=(x_1^1,\ldots,x_{d_1}^1;\ldots;x_1^r,\ldots,x_{d_r}^r)$ with respect
to this decomposition.
Let $C\subset\R^d$ be a convex polyhedral cone with the vertex at the origin $0$.
Assume that $C$ is given by inequalities either of type $x^{i}_{j} \le  a x^{i'}_{j'}$ where $a>0$ and $i\ne i'$ or of type $0\le x^i_j$.
In what follows,
we use the bijective correspondence between facets of $C$ and inequalities, namely, every inequality $x^i_j\le a x^{i'}_{j'}$ defines the facet $H(i,j;i',j')$ given by the equation
$x^i_j=a x^{i'}_{j'}$, and every inequality $0\le x^i_j$ defines the facet $H(0,0;i,j)$ given by the equation $x^i_j=0$.

\begin{figure}
\begin{center}
\includegraphics[width=5cm]{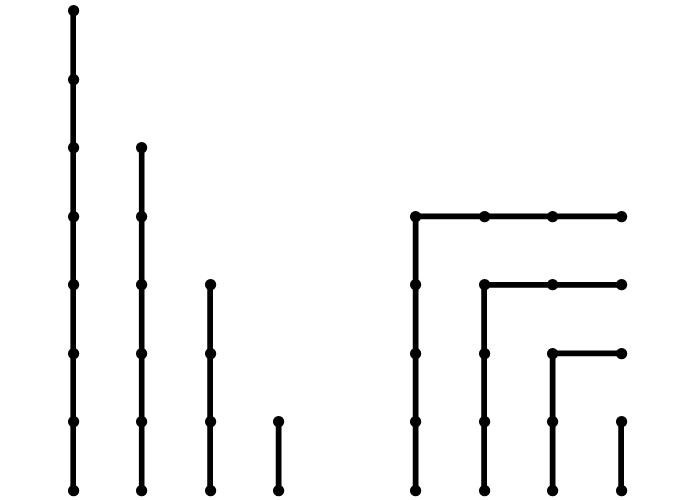}
\hspace{2cm} \label{f.trap}
\includegraphics[width=5cm]{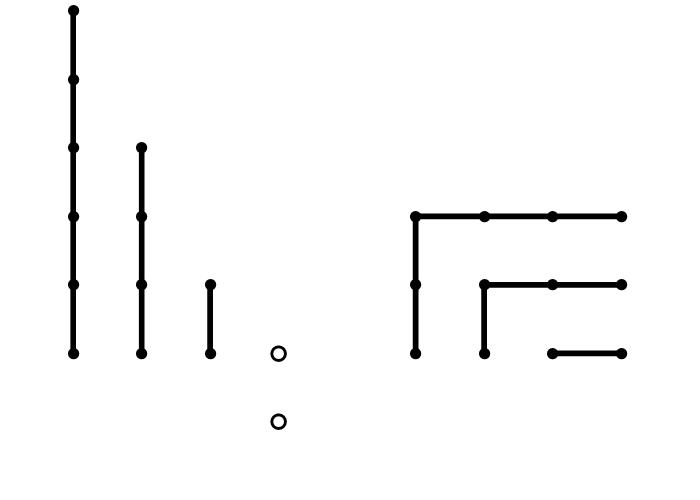}
\caption{}
\end{center}
\end{figure}

In addition, assume that $C$ does not contain any rays parallel to the $x^i_j$-axis unless $j=1$.
Then the geometric mitosis of \cite[Section 5.1]{K16II} can be defined on faces of $C$.
Below we describe the resulting {\em mitosis operations} $M_1$,\ldots, $M_r$ from a
combinatorial viewpoint.

Let $\G$ be a face of the cone $C$ of codimension $\ell$.
The $i$-th {\em mitosis operation} $M_i$ applied to $\G$ will produce a collection $M_i(\Gamma)$ (possibly empty) of faces of $C$.
Choose a minimal subset of facets $H_1$,\ldots, $H_\ell$ of $C$ such that $\G=H_1\cap\ldots \cap H_\ell$.
If none of these facets coincides with $H(p,q;i,d_i)$ for some $p$ and $q$, then $M_i(\G)=\emptyset$.
Otherwise, let $s$ be the smallest number such that the subset $\{H_1,\ldots, H_\ell\}$ contains facets of type
$H(\cdot,\cdot;i,j)$ for all $j=s$, $s+1$, \ldots, $d_i$.
For  brevity, we label these facets by $H^+(i,s)$, $H^+(i,s+1)$, \ldots, $H^+(i,d_i)$.
For every $j=s+1$, $s+2$,\ldots, $d_i$, we now label by $H_+(i,j)$ the facet of type $H(i,j;\cdot,\cdot)$.
If there are two such facets $H(i,j;p,q)$ and $H(i,j;p',q')$, and $x^p_q\le x^{p'}_{q'}$ everywhere on $\G$
then we put $H_+(i,j):=H(i,j;p,q)$.

Let $J_i(\G)\subset\{s, s+1,\ldots,d_i\}$ consist of indices $j$ such that $H_+(i,j)\notin\{H_1,\ldots, H_\ell\}$.
For every $j\in J_i(\G)$, we define an {\em offspring} $\Delta_j\in M_i(\G)$ as the intersection of facets
$$\Delta_j= H_1(j)\cap H_2(j)\cap\ldots\cap H_{\ell-1}(j),$$
where the set $\{H_1(j),\ldots, H_{\ell-1}(j)\}$
is obtained from the set $\{H_1,\ldots, H_\ell\}$ by the following rule.
First, remove the facet $H^+(i,j)$.
Second, for every $k\in J_i(\G)$ such that $k>j$ replace the facet $H^+(i,k)$ by the facet $H_+(i,k)$.
Note that $\dim \Delta_j=\dim \G+1$.
\begin{defin}
The $i$-th {\em mitosis operation} $M_i$ sends $\G$ to
$$M_i(\G)=\{\Delta_j \ | \ j\in J_i(\G)\}.$$
\end{defin}


\subsection{Type $A$: GZ polytopes}
Let $C$ be the vertex cone of the GZ polytope in type $A$ for the vertex
$a=(\l_2,\ldots,\l_n;\l_3,\ldots,\l_n;\ldots;\l_n)$ (see table $(GZ_A)$).
After an affine change of coordinates $x=z-a$ the inequalities that define $C$ can be written as follows:
$$0\le x^1_1; \quad 0\le x^1_2\le x^2_1;\quad \ldots; \quad 0\le x^1_{n-1} \le x^2_{n-2} \le\ldots \le x^{n-1}_1.$$
The cone $C$ has $d=\frac{n(n-1)}{2}$ facets:
$H(0,0;1,i)$ for $1\le i\le (n-1)$ and $H(i-1,j+1;i,j)$ for $2\le i\le (n-1)$, $1\le j\le n-i$.
In particular, we have the following identifications of facets:
$$H(0,0;1,i)=H^+(1,i), \quad H(i-1,j+1;i,j)=H^+(i,j)=H_+(i-1,j+1).$$
It is convenient to encode a face $\G$ of $C$ by an $n\times n$ table ({\em pipe dream}) filled with $+$ as follows.
The table contains $+$ in cell $(i,i+j)$ iff $\G\subset H(i-1,j+1;i,j)$ and $i\ge 2$ or $\G\subset H(0,0;i,j)$ and $i=1$.
In particular, only cells above the main diagonal might have $+$.
In this notation, mitosis operations $M_1$, $M_2$ applied to the vertex $0$ produce the following faces (only cells $(1,2)$, $(1,3)$ and $(2,3)$ of $3\times 3$ tables are shown since the other cells never contain $+$):
$$
\{0\}=
\begin{array}{c}
\begin{array}[t]{|c|}
\hline
+\\
\hline
\end{array}
\begin{array}[t]{|c|}
\hline
+ \\
\hline
+\\
\hline
\end{array}\\
\end{array}
\stackrel{M_1}{\longrightarrow}
\begin{array}{c}
\begin{array}[t]{|c|}
\hline
\ \ \\
\hline
\end{array}
\begin{array}[t]{|c|}
\hline
+ \\
\hline
+\\
\hline
\end{array}\\
\end{array}
\stackrel{M_2}{\longrightarrow}
\begin{array}{c}
\begin{array}[t]{|c|}
\hline
\ \ \\
\hline
\end{array}
\begin{array}[t]{|c|}
\hline
+ \\
\hline
\\
\hline
\end{array}\\
\end{array}
\stackrel{M_1}{\longrightarrow}
\begin{array}{c}
\begin{array}[t]{|c|}
\hline
\ \ \\
\hline
\end{array}
\begin{array}[t]{|c|}
\hline
\ \ \\
\hline
\\
\hline
\end{array}\\
\end{array}=C
$$
$$
\{0\}
\stackrel{M_2}{\longrightarrow}
\begin{array}{c}
\begin{array}[t]{|c|}
\hline
+\\
\hline
\end{array}
\begin{array}[t]{|c|}
\hline
+ \\
\hline
\\
\hline
\end{array}\\
\end{array}
\stackrel{M_1}{\longrightarrow}
\left\{
\begin{array}{c}
\begin{array}[t]{|c|}
\hline
+\\
\hline
\end{array}
\begin{array}[t]{|c|}
\hline
\ \ \\
\hline
\\
\hline
\end{array}\\
\end{array}
\ , \
\begin{array}{c}
\begin{array}[t]{|c|}
\hline
\ \ \\
\hline
\end{array}
\begin{array}[t]{|c|}
\hline
 \\
\hline
+\\
\hline
\end{array}\\
\end{array}
\right\}
\stackrel{M_2}{\longrightarrow}
C$$
(see also Figure 3).

For arbitrary $n$, the mitosis operations $M_1$,\ldots, $M_{n-1}$ encoded by tables coincide with Knutson--Miller mitosis on pipe dreams \cite{KnM}
after reflecting tables in a vertical line.
Instead of the vertex cone $C$ we could take the GZ polytope in type $A$ and consider mitosis operations on faces that contain the vertex $a$ (so called {\em Kogan faces}).
Geometric meaning of the resulting collections of faces is described in \cite[Theorem 5.1, Corollary 5.3]{KST}.
In particular, the following analog of Kushnirenko's theorem holds.

Recall that Schubert subvarieties $X_w$ are labeled by the elements of the Weyl group of $G$, namely,
$X_w$ is the closure of the $B$-orbit $BwB/B$, where $w$ is an element of the Weyl group of $G$.
The Weyl group of $G=SL_n(\C)$ is the symmetric group $S_n$.
By $s_1$,\ldots, $s_{n-1}$ we denote the elementary transpositions.
\begin{thm}{\cite[Theorem 5.4]{KST}}
Let $X_w\subset SL_n(\C)/B$ be the Schubert subvariety corresponding to a permutation $w\in S_n$.
Let $w=s_{j_1}\ldots s_{j_\ell}$ be a reduced decomposition of a permutation $w\in S_n$ such that
$(j_1,\ldots,j_\ell)$ is a subword of $(1;2,1;3,2,1;\ldots;{n-1},\ldots,1)$.
Let $\Sch_w\subset GZ_\l$ be the set of all faces produced from the vertex $a\in GZ_\l$ by applying successively the operations $M_{n-j_\ell}$,\ldots, $M_{n-j_1}$:
$$\Sch_w=\ M_{n-j_1}\cdots M_{n-j_\ell}(a).$$
Then
$$\deg_\l(X_w)=\ell! \sum_{\G\in \Sch_w}\Volume(\G).$$
\end{thm}
This implies that the Schubert cycle $[X_w]$ (that is, the cohomology class of $X_w$ in $H^*(SL_n(\C)/B,\Z)$) in the polytope ring $R_{GZ}\simeq H^*(SL_n(\C)/B,\Z)$ is represented by the sum of faces in $\Sch_w$.

\begin{figure}
\begin{center}
\includegraphics[width=2.5cm]{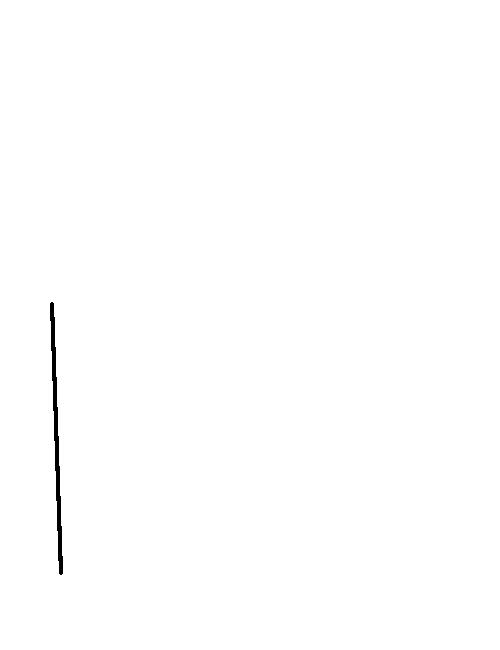}
\label{f.GZ_faces}
\includegraphics[width=2.5cm]{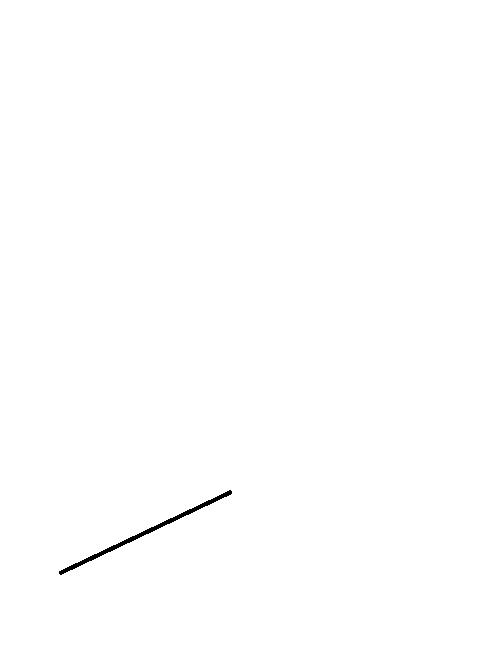}
\includegraphics[width=2.5cm]{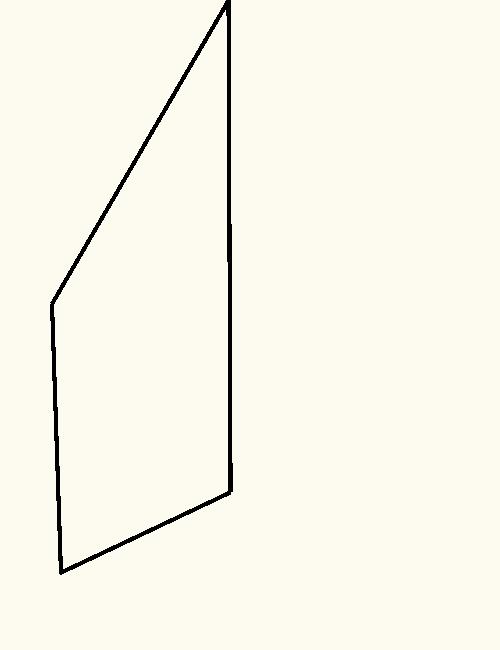}
\includegraphics[width=2.5cm]{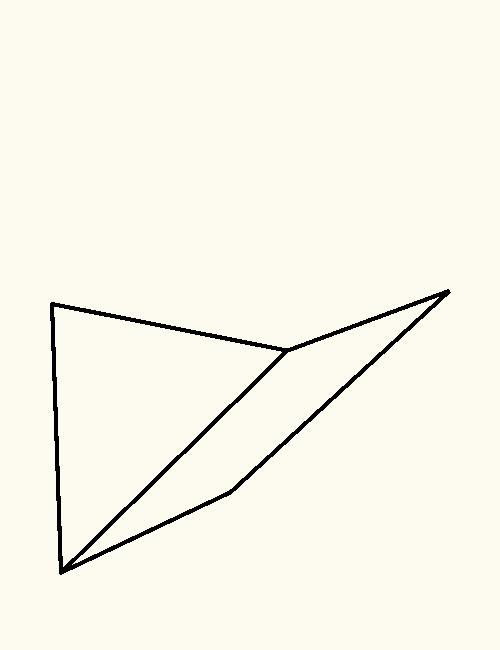}
\includegraphics[width=2.5cm]{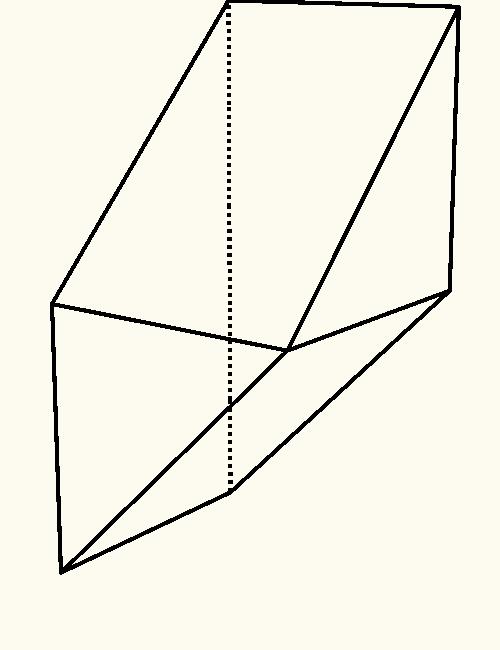}
\caption{Faces $M_1(a)$, $M_2(a)$, $M_2M_1(a)$, $M_1M_2(a)$, $M_1M_2M_1(a)$ of the GZ polytope in type $A$, $n=3$}
\end{center}
\end{figure}

\begin{example} For $n=3$, we have $[X_{s_1s_2}]=M_2M_1(a)$ and $[X_{s_2s_1}]=M_1M_2(a)$.
Since the faces in these two presentations are transverse and their intersection consists of two edges $M_1(a)$ and $M_2(a)$ we get the identity:
$[X_{s_1s_2}]\cdot[X_{s_2s_1}]=[X_{s_1}]+[X_{s_2}]$ (see Figure 3).
\end{example}

\subsection{Type $C_{2-3}$: DDO polytopes} \label{ss.C_2}
In \cite[Example 2.9]{K16I}, the following family of DDO polytopes in $\R^4=\R^2\oplus\R^2$ is considered:
$$0\le x_1^1\le \l_1, \quad x^2_1\le x_1^1+\l_2, \quad x^1_2\le 2x^2_1,$$
$$x^1_2\le x^2_1+\l_2,\quad 0\le x^2_2\le \l_2, \quad x^2_2\le \frac{x^1_2}{2}$$
(these polytopes can also be realized as Newton--Okounkov polytopes of the isotropic flag variety $Sp_4/B$ \cite[Proposition 4.1]{K16I})
The vertex cone $C$ of the vertex $0$ is given by 4 homogeneous inequalities:
$0\le x^1_1$, $0\le 2x^2_2\le x^1_2\le 2x^2_1$.
It is convenient to encode a face $\G$ of $C$ by a $(2n-1)\times n$ table ({\em skew pipe dream}) for $n=2$ filled with $+$ as follows (see Section \ref{ss.C} for the general definition of skew pipe dreams).
The table contains $+$ in cell $(3-i,i)$ (for $i=1$, $2$) iff $\G\subset H(0,0;i,i)$,
$+$ in cell $(2,2)$ iff $\G\subset H(2,2;1,2)$ and $+$ in cell $(3,2)$ iff $\G\subset H(1,2;2,1)$.
There are two mitosis operations $M_1$ and $M_2$.
$$\{0\}=
\begin{array}{|c|}
\hline
+\\
\hline
\end{array}\begin{array}{|c|}
\hline
+\\
\hline
+ \\
\hline
+ \\
\hline
\end{array}
\stackrel{M_1}{\longrightarrow}
\begin{array}{|c|}
\hline
\ \ \\
\hline
\end{array}\begin{array}{|c|}
\hline
+\\
\hline
+ \\
\hline
+ \\
\hline
\end{array}
\stackrel{M_2}{\longrightarrow}
\begin{array}{|c|}
\hline
\ \ \\
\hline
\end{array}\begin{array}{|c|}
\hline
+\\
\hline
+ \\
\hline
 \\
\hline
\end{array}
\stackrel{M_1}{\longrightarrow}
\begin{array}{|c|}
\hline
\ \ \\
\hline
\end{array}\begin{array}{|c|}
\hline
+\\
\hline
 \\
\hline
 \\
\hline
\end{array}
\stackrel{M_2}{\longrightarrow}
\begin{array}{|c|}
\hline
\ \ \\
\hline
\end{array}\begin{array}{|c|}
\hline
\\
\hline
\ \  \\
\hline
\\
\hline
\end{array}=C$$
$$\{0\}
\stackrel{M_2}{\longrightarrow}
\begin{array}{|c|}
\hline
+\\
\hline
\end{array}\begin{array}{|c|}
\hline
+\\
\hline
+ \\
\hline
 \\
\hline
\end{array}
\stackrel{M_1}{\longrightarrow}
\left\{\begin{array}{|c|}
\hline
\ \ \\
\hline
\end{array}\begin{array}{|c|}
\hline
+\\
\hline
 \\
\hline
+ \\
\hline
\end{array}\ , \
\begin{array}{|c|}
\hline
+ \\
\hline
\end{array}\begin{array}{|c|}
\hline
+\\
\hline
\\
\hline
 \\
\hline
\end{array}
\right\}
\stackrel{M_2}{\longrightarrow}
$$
$$
\stackrel{M_2}{\longrightarrow}\left\{\begin{array}{|c|}
\hline
\ \ \\
\hline
\end{array}\begin{array}{|c|}
\hline
\\
\hline
+ \\
\hline
 \\
\hline
\end{array} \ , \
\begin{array}{|c|}
\hline
\ \ \\
\hline
\end{array}\begin{array}{|c|}
\hline
\\
\hline
 \\
\hline
+ \\
\hline
\end{array}
\ , \
\begin{array}{|c|}
\hline
+\\
\hline
\end{array}\begin{array}{|c|}
\hline
\\
\hline
\ \  \\
\hline
\\
\hline
\end{array}\right\}
\stackrel{M_1}{\longrightarrow}C
$$
The Weyl group of $G=Sp_4(\C)$ is the dihedral group $D_4$.
By $s_1$, $s_2$ we denote simple reflections that generate $D_4$ so that $s_2$ corresponds to the longer root.
By \cite[Corollary 3.6]{K16I} we have
\begin{prop}\label{p.DDO_C}
Let $X_w\subset Sp_4(\C)/B$ be the Schubert subvariety corresponding to a permutation $w\in D_4$.
Let $w=s_{j_1}\ldots s_{j_\ell}$ be a reduced decomposition of a permutation $w\in D_4$ such that
$(j_1,\ldots,j_\ell)$ is a subword of $(1,2;1,2)$.
Let $\Sch_w\subset GZ_\l$ be the set of all faces produced from the vertex $a\in GZ_\l$ by applying successively the operations $M_{j_\ell}$,\ldots, $M_{j_1}$:
$$\Sch_w=\ M_{j_1}\cdots M_{j_\ell}(0).$$
Then
$$\deg_\l(X_w)=\ell! \sum_{\G\in \Sch_w}\Volume(\G).$$
\end{prop}

This example can be extended to DDO polytopes for $Sp_{2n}$.
For $n=3$ and the DDO polytope for $(s_3s_2s_1)^3$ (where $s_3$ is the simple reflection with respect
to the longer root) this was done in \cite{P15}.
The corresponding family of DDO polytopes in $\R^9=\R^3\oplus\R^3\oplus\R^3$ is given by inequalities:
$$0\le x_1^1\le \l_1; \quad x_1^2\le \l_2+x_1^1; \quad x_1^3\le \l_3+x_1^2;$$
$$0\le x_2^1\le\min\{x_1^2,\l_2\}; \quad x_2^2\le\min\{\l_3+x_2^1+x_1^3,2x_1^3\};$$
$$x_2^3\le \min\{x_2^1+\l_3,\frac12x_2^2\};\quad
x_3^1\le\min\{x_2^2,x_1^3+\l_3,\l_3+x_2^2-x_2^3\};$$
$$x_3^2\le\min\{x_3^1,x_2^3+\l_3,2x_2^3\}; \quad 0\le x_3^3\le\min\{\frac12x_3^2,\l_3\}.$$
In particular, the vertex cone at $0$ is not simplicial.
It is defined by $10$ inequalities:
$$0\le x^1_1; \quad 0\le x^1_2\le x^2_1; \quad 0\le x^3_3\le\frac12 x^2_3\le \frac12 x^1_3 \le\frac12 x^2_2\le x^3_1;$$
$$\frac12x^2_3\le x^3_2\le\frac12 x^2_2.$$
An analog of Proposition \ref{p.DDO_C} follows easily from \cite[Corollary 3.6]{K16I}.
However, combinatorics of mitosis becomes more involved as analogs of pipe dreams in this case have a loop.

Recently, Fujita identified DDO polytopes with certain Nakashima--Zelevinsky polyhedral realizations of crystal bases \cite[Theorem 4.1]{F}.
In particular, there are explicit inequalities for these polytopes in types $A$, $B$, $C$, $D$ and $G_2$ \cite[Example 4.3]{F}.
In type $A$, they coincide with the GZ polytope and in type $C_{2-3}$ with the polytopes described in this section.
It would be interesting to apply geometric mitosis to these polytopes in the other cases.

\subsection{Type $C$: GZ polytopes}\label{ss.C}
The combinatorics of $C_2$ example from Section \ref{ss.C_2} can be extended to $C_n$ in a different way by using the string cone $C$ for the reduced decomposition $\overline{w_0}=(s_n s_{n-1}\ldots s_2 s_1 s_2\ldots s_{n-1}s_n)\ldots (s_2s_1s_2)(s_1)$ of the longest element in the Weyl group (here $s_1$ corresponds to the longer root).
The corresponding string polytope coincides with the symplectic GZ polytope after a unimodular change of coordinates \cite[Section 6]{L}.
The cone $C$ is simplicial and is given by $d=n^2$ inequalities:
$$0\le x^i_2\le x^{i-1}_4\le x^{i-2}_6\le\ldots\le x^2_{2i-2}\le x^1_{i}\le x^2_{2i-3}\le\ldots\le x^{i-2}_5\le x^{i-1}_3\le x^i_1$$
for all $i=1$,\ldots, $n$.
We define {\em symplectic mitosis} as the geometric mitosis associated with the cone $C$.
Combinatorics of the symplectic mitosis is quite simple and described in detail in \cite[Section 5.2]{K16I} using {\em skew pipe dreams}.
However, arguments of \cite[Corollary 3.6]{K16I} do not directly yield presentations for Schubert cycles since the symplectic GZ polytope does not satisfy the necessary conditions.
Still computations for $n=2,3$ suggest that the collections of faces of the symplectic GZ polytope obtained using symplectic mitosis do represent the corresponding Schubert cycles in the polytope ring $R_{SGZ}$.
Below we describe a bijection between faces of $C$ and faces of $SGZ_\l$ that we used.

Let $v$ be the vertex of $SGZ_\l$ given by equations $\l_s=x^i_j=y^k_l$ for all triples $\l_s$, $x^i_j$ and $y^k_l$ such that
$s=i+j-1=k+l$.
We now define a bijection between those facets of $P_\l$ that contain $v$ and {\em skew pipe dreams}
of size $n$ with exactly one $+$.
Recall that a {\em skew pipe dream} of size $n$ is a $(2n-1)\times{n}$ table whose cells are either empty or filled with $+$.
Only cells $(i,j)$ with $n-j< i < n+j$ are allowed to have $+$ (see \cite[Section 5.2]{K16I} for more details on skew pipe dreams).
Put $y^0_i:=\l_i$ for $i=1$,\ldots, $n$.
The facet given by equation  $x^i_j=y^{i-1}_j$ corresponds to the skew pipe dream with $+$ in cell $(i+j-1,n-i+1)$.
The facet given by equation $y^i_j=x^i_{j+1}$ corresponds to the skew pipe dream with $+$ in cell $(2n-i-j+1,n-i+1)$.
In what follows, we denote by $H_{(k,l)}$ the facet whose skew pipe dream under this correspondence contains $+$
in cell $(k,l)$.

This correspondence between facets and skew pipe dreams with a single $+$ extends to all faces of the symplectic GZ polytope that contain the vertex $v$.
Namely, the face $H_{k_1,l_1}\cap\ldots\cap H_{k_s,l_s}$ obtained as the intersection of $s$ facets corresponds to the skew pipe dream that has $+$ precisely in cells $(k_1,l_1)$,\ldots, $(k_s,l_s)$.
In particular, the vertex $v$ corresponds to the skew pipe dream $D_0$ that has $+$ in all (fillable) cells.
In what follows, we denote by $F_D$ the face corresponding to a skew pipe dream $D$.

We now formulate a conjecture.
Let $w$ be an element of the Weyl group of $G=Sp_{2n}$.
Choose  a reduced decomposition $w=s_{j_1}\ldots s_{j_\ell}$
such that it is a subword of $(s_n s_{n-1}\ldots s_2 s_1 s_2\ldots s_{n-1}s_n)\ldots (s_2s_1s_2)(s_1)$.

\begin{conj} Define the set $\Sch_w$ of faces of the symplectic GZ polytope as follows:
$$\Sch_w=\{F_D \ |\ D\in M_{n+1-j_\ell}\cdots M_{n+1-j_1}(D_0)\}$$
where $M_i$ denotes the $i$-th symplectic mitosis operation.
Then
$$\deg_\l(X_w)=\ell!\sum_{F\in \Sch_w} \Volume(F).$$
\end{conj}

This conjecture is verified in the case $n=2$ and for certain $w$ in the case $n=3$ \cite{P17}.
Note that the bijection between faces of $SGZ_\l$ that contain the vertex $v$ and faces of the string cone $C$ does not come from the unimodular change of coordinates that identifies the string polytope and the symplectic GZ polytope.
There are might be piecewise linear transformations (such as the ones used in \cite[Section 5.2]{K15}) that yield scissor congruence of unions of faces of $SGZ_\l$ and faces of another  polytope for which geometric meaning of symplectic mitosis is more transparent.

\subsection{Type $B$: GZ polytopes}
Note that the Weyl groups of $Sp_{2n}(\C)$ and $SO_{2n+1}(\C)$ are the same.
Since the GZ polytopes for both groups differ only by lattices symplectic mitosis is also a natural tool for finding presentations of Schubert cycles by faces of $OGZ_\l$ in type $B$.
However, coefficients will be rational rather than integer (with powers of 2 in denominator) because the torsion index of $SO_{2n+1}(\C)$ is a power of 2.
Note also that the volumes of faces of both $SGZ_\l$ and $OGZ_\l$ should be computed with respect to their lattices.
The difference is already visible in the case $n=2$ (see Example \ref{e.B_2}).

\end{document}